\documentclass[journal]{IEEEtran}
\ifCLASSINFOpdf
\else
\fi
\usepackage{amsmath,amsfonts,amsthm}

\newtheorem{theorem}{Theorem}
\newtheorem{lemma}{Lemma}
\theoremstyle{remark}
\newtheorem{rem}{Remark}

\usepackage[normalem]{ulem}
\newcommand\redout{\bgroup\markoverwith
{\textcolor{red}{\rule[.5ex]{2pt}{0.4pt}}}\ULon}
\usepackage{xcolor}

\usepackage{algorithm,algorithmic}
\begin{document}
%
\title{Distributed Banach-Picard Iteration for Locally Contractive Maps}
%
%
%

\author{Francisco~Andrade,
        M\'{a}rio~A.~T.~Figueiredo,~\IEEEmembership{Fellow,~IEEE},
        and Jo\~{a}o~Xavier
\thanks{The authors are with the Instituto Superior T\'{e}cnico, Universidade de Lisboa, Lisbon, Portugal.} %
\thanks{F.~Andrade (francisco.andrade@tecnico.ulisboa.pt) and M.~Figueiredo (mario.figueiredo@tecnico.ulisboa.pt) are also with the Instituto de Telecomunicaç\~{o}es, Lisbon, Portugal.}
\thanks{J.~Xavier (jxavier@isr.tecnico.ulisboa.pt) is  also with the Laboratory for Robotics and Engineering Systems, Institute for Systems and Robotics, Lisbon, Portugal.}}

%
%

\markboth{Submitted}%
{Shell \MakeLowercase{\textit{et al.}}: Bare Demo of IEEEtran.cls for IEEE Journals}
%



\maketitle

\begin{abstract}
The Banach-Picard iteration is widely used to find fixed points of locally contractive (LC) maps. This paper extends the Banach-Picard iteration to distributed settings; specifically, we assume the map of which the fixed point is sought to be the average of individual (not necessarily LC) maps held by a set of agents linked by a communication network. An additional difficulty is that the LC map is not assumed to come from an underlying optimization problem, which prevents exploiting strong global properties such as convexity or Lipschitzianity. Yet, we propose a distributed algorithm and prove its convergence, in fact showing that it maintains the linear rate of the standard Banach-Picard iteration for the average LC map. As another contribution, our proof imports tools from perturbation theory of linear operators, which, to the best of our knowledge, had not been used before in the theory of distributed computation.

\end{abstract}


%
\IEEEpeerreviewmaketitle

\section{Introduction}

\IEEEPARstart{I}{nferring}  a desired quantity $x^\star$ from data, can often be naturally expressed as a fixed point equation $H(x^\star)=x^\star$,  the map $H$  relating the data to the desired quantity. Although in some rare  cases this fixed point equation can be solved in closed-form, more often than not, $x^\star$ has to be numerically approximated using, \textit{e.g.}, the so-called  \emph{Banach-Picard iteration}:  
\begin{align}
x^{k+1}=H(x^k). \label{Picard}
\end{align}
For a recent comprehensive review of the fixed-point strategy to inference problems, see \cite{Compbettes2020}.

If the entire data, thus the operator $H$,  is available to some agent, that agent can  perform the Banach-Picard iteration. In contrast, in the so-called \textit{distributed} scenario, the data is acquired by spatially dispersed agents who only have access to local data. In such distributed setups, no single agent possesses the full data set, hence no single agent can compute the map $H$. Instead, each agent holds a local portion of the data and can communicate only with a subset of the other agents (its neighbours). Nevertheless, the goal remains that of finding a fixed point of $H$, under the constraints of this  distributed configuration: each agent can only engage in private/local computation and in communication with its neighbours. 

\subsection{Problem Statement} \label{Problem Statmente}
We consider a network of $N$ agents, where the interconnection structure is represented by an undirected and connected graph: the nodes correspond to the agents and an edge between two agents indicates  they can communicate (are neighbours). Each agent $n \in \{1,...,N\}$ holds an operator $H_n:\mathbb{R}^d \to \mathbb{R}^d$, and the goal is to compute a fixed point of the average operator 
\begin{align}
H=\frac{1}{N}\sum_{n=1}^N H_n. \label{DistributedMap}
\end{align}
Crucially, each agent $n$ is restricted to performing computations involving $H_n$ and communicating with its neighbours. 

Our only assumption about $H$ (which may not hold for each $H_n$)  is that it has an \emph{attractor}, {\it i.e.}, a fixed point $x^\star$, such that $H$ is continuously differentiable in a neighbourhood of $x^\star$ and 
\begin{align} \label{JacobianCondition}
\rho\big(\mathbf{J}_H(x^\star)\big)<1,
\end{align}
where $\rho\big(\mathbf{J}_H(x^\star)\big)$ is the spectral radius of the Jacobian of $H$ at $x^\star$. In other words, $H$ is \textit{locally contractive} (LC) in a neighbourhood of $x^\star$.

Application of Ostrowski's Theorem (see Appendix \ref{Auxiliary results Section}) to the Banach-Picard iteration \eqref{Picard} ensures that if the sequence $x^k$ gets sufficiently close to $x^\star$, then it not only converges to $x^\star$, but does so at a linear rate; \textit{i.e.}, it exhibits local linear convergence.

\subsection{Contributions}
Our main contribution is a distributed algorithm that, as the corresponding centralized Banach-Picard iteration \eqref{Picard}, has local linear convergence at attractors of $H$. 	

Although assuming a relatively weak set of conditions---essentially only local linear convergence of the centralized Banach-Picard iteration---and no global structure (\textit{e.g.}, Lipschitzianity or coercivity), we propose a distributed algorithm and prove that it inherits the local linear convergence of its centralized counterpart. 

Even though the assumptions are rather weak, they nevertheless suffice to encapsulate relevant algorithms, namely some instances of the \textit{expectation maximization} (EM, \cite{McLachlan2007}, \cite{pereira2018parameter}) algorithm and the one proposed in \cite{sanger1989optimal} for \textit{principal component analysis} (PCA). A forthcoming paper is devoted to the application of the algorithmic framework herein proposed  to obtain distributed versions of those algorithms, with  local linear convergence guarantees  \cite{FranciscodistributedEMPCA}.

As an additional contribution, we mention the proof technique, which, as far as we know, strongly departs from the standard proof techniques used in distributed optimization. Specifically, we employ tools from perturbation theory of linear operators \cite{kato2013perturbation}, which, to the best of our knowledge, had not been used before in the context of distributed computation. Arguably, there are proof techniques that resemble a ``perturbative
argument'' on the eigenvalues of a matrix (e.g., Proposition 2.8 in
\cite{blatt2007convergent}). However, those techniques bypass the subtle issue of the
differentiability of the eigenvalues, simply using the formula for the
derivative of the determinant. In contrast, the theorem from
perturbation theory (PT) of linear operators that we use
simultaneously handles the differentiability issue and simplifies the
computation of the derivative.  To the best of our knowledge, the full
power of PT had never been used in the context of distributed
optimization.
\subsection{Remarks}

Our setup departs from standard ones in two main aspects. First, it encompasses problems that are not naturally expressed as optimization problems. This last notion should be understood with a grain of salt, since a fixed point of $H$ minimizes $\|H(x)-x\|^2$;  however, in many cases, there is a more ``natural'' objective function  than this one. For example, if the Jacobian $\mathbf{J}_H$ is symmetric in an open, convex set, then there exists a function $f$ such that $H=\nabla f$ \cite[Theorem 1.3.1]{facchinei2007finite}, and the Banach-Picard iteration can be seen as method to find a stationary point of $ f(x) -\frac{1}{2}\|x\|^2$.   

Second, the notion of  attractor herein used is purely local. By assuming \eqref{JacobianCondition}, we consider only  local guarantees. Many optimization problems benefit from global properties, such as Lipschitzianity or strong convexity. Such properties, however, are absent from many relevant algorithms, such as EM, for which only local guarantees can be given.

\subsection{Related Work}
In this section, we review relevant related work in distributed computation, highlighting how our contributions differ from that other work.

A setup that closely resembles ours is considered in \cite{gang2019fast} and \cite{pereira2018parameter}; in fact, the problems therein addressed are, respectively, distributed PCA and distributed EM. As shown in the upcoming paper \cite{FranciscodistributedEMPCA}, our setup encapsulates the problems addressed in \cite{gang2019fast} and \cite{pereira2018parameter}. However, the algorithm in  \cite{pereira2018parameter} uses a diminishing step-size, which, unlike our algorithm, results in a sacrifice of the convergence rate   of the centralized EM. The algorithm in \cite{gang2019fast} is recovered by using our approach to build a distributed version of the algorithm in  \cite{sanger1989optimal}. Moreover, our work has at least two advantages over that of  \cite{gang2019fast}: we provide a proof of local linear convergence (which \cite{gang2019fast} does not) and our setup is not restricted to the algorithm in \cite{sanger1989optimal}.

The works \cite{fullmer2018distributed}, \cite{li2020linearly}, and \cite{li2020distributed} share  a similarity with ours by addressing the distributed computation of fixed points. However, the setups therein considered have much more structure than ours: Lipschitzianity and quasi-nonexpansivness \cite{li2020linearly}, non-expansiveness \cite{li2020distributed}, and paracontractiveness \cite{fullmer2018distributed}. Those are global properties which are absent in algorithms such as EM or the algorithm in \cite{sanger1989optimal} for PCA.

A  large body of work on distributed optimization has been produced in the last decade; see \cite{nedic2009distributed}, \cite{Jakovetic2014}, \cite{Jakovetic2015},  \cite{shi2015extra}, \cite{xu2015augmented}, \cite{giannakis2016decentralized}, \cite{qu2017harnessing}, \cite{nedic2017achieving}, \cite{xu2017adaptive},  \cite{jakovetic2018unification}, \cite{alghunaim2020linear},\cite{uribe2020dual},  \cite{qu2019accelerated}, \cite{mansoori2019general}, \cite{jakovetic2020primal},  and \cite{fallah2019robust}  for  convex optimization, where the last reference considers a stochastic variant, and \cite{di2016next}, \cite{tatarenko2017non}, \cite{vlaski2021distributed}, as  examples of distributed non-convex optimization. All the algorithms in those works can definitely be seen as distributed algorithms for finding fixed points. However, as their setups stem from optimization, they further assume conditions such as coercivity, Lipschitzianity, or strong convexity.  In our work, none of these properties are assumed, and  only a basic local assumption is made.

\subsection{Organization of the Paper}
This introductory section will conclude with a brief paragraph describing the adopted notation. Section \ref{KKT section} describes the proposed distributed algorithm. The theorems guaranteeing local linear convergence  are presented and proved in Section \ref{Convergence}. Section \ref{Intuition} provides some intuitive insight into the proposed algorithm. The paper concludes in Section \ref{Conclusions}, which also points at ongoing and future work. The appendix contains three theorems that are instrumental in this work.

\subsection*{Notation}
The set of real $n$ dimensional vectors with positive components is denoted by $\mathbb{R}^n_{>0}$.
Matrices and vectors are denoted by upper and lower case letters, respectively. The spectral radius of a matrix $M$ is denoted by $\rho(M)$. Given a map $H$ and a function $f$, $\mathbf{J}_H(x)$ and $\mathbf{J}^2_f(x)$ denote, respectively, the Jacobian of $H$ and the Hessian of $f$, at $x$. Given a vector $v$, $v_s$ denotes its $s$th component; given a matrix $M$, $M_{st}$ denotes the element on the $s$th line and $t$th column; $M^T$ is its transpose. 
The $d$-dimensional identity matrix is denoted by $I_d$, and $\mathbf{1}_d$ is the $d$-dimensional vector of ones. The Kronecker product is denoted by $\otimes$. The letter $i$ is denotes the imaginary unit ($i^2=-1$). The fixed point set of a map $H$ is denoted by $\mbox{Fix}(H)$. Given a matrix $L$, its Moore-Penrose (pseudo)inverse is denoted as $L^+$. The gradient of a function $f$ with respect to vector $w$ is denoted by $\nabla_w f$. If a matrix $M$ is negative (positive) definite, this is denoted by $M \prec 0$ ($M \succ 0$). Whenever convenient, we will denote a vector with two stacked blocks, $[v^T, u^T]^T$, simply as $(v,u)$.

\section{The Distributed Algorithm} \label{KKT section}
Let $R$ be the map on $\mathbb{R}^{dN}$ defined, for $z=[z_1^T,\ldots,z_N^T]^T$ with $z_j \in \mathbb{R}^d$, by 
$$
R(z)=\Big[ \big(H_1(z_1)-z_1\big)^T,\ldots,\big(H_N(z_N)-z_N\big)^T\Big]^T,
$$ and let $W=\tilde{W}\otimes I_d$, where $\tilde{W}$ is the Metropolis weight matrix associated to the agents' communication graph  \cite{xiao2004fast}. The distributed algorithm we propose is presented in Algorithm 1, where $\alpha \in \mathbb{R}_{>0}$.

\begin{algorithm}
\caption{Distributed Banach-Picard Iteration}
\begin{algorithmic}[1]
  \small
  \STATE Initialization:
  \begin{equation*}
\begin{aligned}
z^{0} &\in \mathbb{R}^{dN},\\
z^1 &= Wz^0+\alpha R(z^0),\\
\end{aligned}
\end{equation*}
\STATE Update:
\begin{equation*}
z^{k+2}=(I+W)z^{k+1}-\frac{I+W}{2}z^k+\alpha\big(R(z^{k+1})-R(z^k)\big).
\end{equation*} 
\end{algorithmic}
\end{algorithm}

This algorithm, inspired by EXTRA (see \cite{shi2015extra} and section \ref{Intuition}) below), is a particular instance of the parametric family of distributed algorithms given by
\begin{equation} \label{ParametricAlgorithm}
\begin{aligned}
z^{0} &\in \mathbb{R}^{dN},\\
z^1 &= z^0+\alpha R(z^0)-\eta L z^0,\\
z^{k+2}&=(2I-\eta L)z^{k+1}-(I+\beta
^2L-\eta L)z^k\\
&+\alpha\big(R(z^{k+1})-R(z^k)\big),
\end{aligned}
\end{equation}
together with the following assumptions on the $dN \times dN$ matrix $L$:
\begin{description}
	\item[a)] $L$ is symmetric and positive semidefinite; 
	\item[b)] $\rho(L)<2$;
	\item[c)] $\ker(L)=\lbrace z \in \mathbb{R}^{dN}: z_1=\cdots=z_N\rbrace$; 
	\item[d)] $L=\tilde{L}\otimes I_d$, where $\tilde{L}$ is $N \times N$ and has the property that $\tilde{L}_{st}=0$ if agent $s$ does not communicate with agent $t$ (thus establishing the compatibility of the algorithm with the network structure).
\end{description}

In fact, Algorithm 1 is recovered from \eqref{ParametricAlgorithm} by making the choices $\eta=1$, $\beta^2=\frac{1}{2}$, and $L=I-W$ (which satisfies assumptions a) -- d) above).

In the next section we establish the local linear convergence of \eqref{ParametricAlgorithm} and, from  remark 1 (see \ref{SectionTheoremStatements}), the instance yielding Algorithm 1.

\section{Convergence Analysis} \label{Convergence}
\subsection{Auxiliary maps} \label{Fmaps}
To study the convergence of \eqref{ParametricAlgorithm} and, consequently of Algorithm 1, we introduce two auxiliary maps, whose connection with Algorithm 1 is explained in Section \ref{SectionOnconnection}. Let $\tilde{U}$ be a matrix with columns forming an orthonormal basis of $\mbox{range}(L)$ and let $F:\mathbb{R}^{dN}\! \times \mathbb{R}^{dN}\! \to \mathbb{R}^{dN}\! \times \mathbb{R}^{dN}$ and $\tilde{F}:\mathbb{R}^{dN}\! \times \mathbb{R}^{d(N-1)}\! \to \mathbb{R}^{dN}\! \times \mathbb{R}^{d(N-1)}$ be the maps defined as, respectively, 
\begin{align}
F\Big(\begin{bmatrix}
z\\
w
\end{bmatrix}\Big) & =
\begin{bmatrix}
z+\alpha R(z)+\beta L^\frac{1}{2}w-\eta Lz\\
w-\beta L^\frac{1}{2}z
\end{bmatrix}  \label{Fmap}
\end{align}
and
\begin{align}
\tilde{F}\Big(\begin{bmatrix}
z\\
\tilde{w}
\end{bmatrix}\Big) &=
\begin{bmatrix}
z+\alpha R(z)+\beta L^\frac{1}{2}\tilde{U}\tilde{w}-\eta Lz\\
\tilde{w}-\beta \tilde{U}^TL^\frac{1}{2}z
\end{bmatrix}.
\end{align}
In section \ref{SectionOnconnection}, it is shown that the elimination of the second variable in the Banach-Picard iteration
\begin{align}
\begin{bmatrix}
z^{k+1}\\
w^{k+1}
\end{bmatrix}=F\Big(\begin{bmatrix}
z^k\\
w^k
\end{bmatrix}\Big)\ \label{Picard2}
\end{align}
recovers algorithm \eqref{ParametricAlgorithm}. The main convergence theorem of this paper (Theorem \ref{Corolary}) can be informally  stated as follows: it is possible to choose $\alpha>0, \beta>0$, and $\eta>0$, such that if $z^k$ in \eqref{Picard2} becomes sufficiently close to $\mathbf{1}_N\otimes x^\star$, then $z^k$ converges to $\mathbf{1}_N\otimes x^\star$ with at least linear rate. The rest of the section is devoted to making this notion precise.

Observe that if $(z^\star,w^\star)$ is a fixed point of $F$, then $z^\star=\mathbf{1}_N\otimes x^\star$, where $x^\star$ is a fixed point of $H$. Conversely, if $x^\star$ is a fixed point of $H$, then, there exists $w^\star$ such that $(\mathbf{1}_N \otimes x^\star, w^\star)$ is a fixed point of $F$. However, note as well that if $(z^\star,w^\star)$ is a fixed point of $F$, then any point in $(z^\star,w^\star+\ker (L))$ is also a fixed point of $F$. As a consequence, for each fixed point of $F$, there is an affine subspace  of fixed points that contains it. A simple corollary of this observation, together with Ostrowski's Theorem (see Theorem \ref{Ostrowski's THeorem} in Appendix \ref{Auxiliary results Section}), is that   $(z^\star,w^\star)$ is not an attractor of $F$ (no fixed point of $F$ is isolated). The way around this issue is to ``identify''\footnote{The meaning of ``identify'' should be understood as an equivalence relation. By considering $\mathbb{R}^{dN}\times \mathbb{R}^{dN}$ modulo $\lbrace 0 \rbrace \times \ker(L)$, we can treat a set of fixed points of the form $\big(z^\star,w^\star+\ker(L)\big)$ as  an equivalence class. This notion of  quotient in case of vector spaces can be easily ``done in coordinates'' by introducing the matrix $\tilde{U}$. In fact, $\mathbb{R}^{dN}$ modulo $\ker(L)$ is isomorphic to $\mbox{range}(L)$.} all the fixed points in $(z^\star, w^\star+\ker (L))$; this is the role of the map $\tilde{F}$.

Note that $\tilde{F}$ identifies the fixed points of $F$ of the form $(z^\star,w^\star+\ker (L))$ and  there is a one-to-one correspondence between the $\mbox{Fix}(H)$ and $\mbox{Fix}(\tilde{F})$. In fact, note that the map $\psi$ on $\mbox{Fix}(H)$ defined by
\begin{align*}
\psi(x^\star)=(\mathbf{1}_N\otimes x^\star,\tilde{w}^\star),
\end{align*}
where $\tilde{w}^\star$ is the unique point satisfying $$\frac{\alpha}{\beta}(L^\frac{1}{2})^+R(\mathbf{1}_N\otimes x^\star)=-\tilde{U}\tilde{w}^\star,$$ is a bijection between $\mbox{Fix}(H)$ and $\mbox{Fix}(\tilde{F})$.

Section \ref{SectionTheoremStatements} shows that $\psi$ preserves condition \eqref{JacobianCondition}, in the sense that if $x^\star$ is an attractor of $H$, then $\psi(x^\star)$ is an attractor of $\tilde{F}$.

\subsection{Connection between $F$ and  \eqref{ParametricAlgorithm}} \label{SectionOnconnection}
The presence of matrix $L^{\frac{1}{2}}$ prevents the map $F$ defined in \eqref{Fmap}, and which defines the iteration \eqref{Picard2}, from having a distributed implementation (whereas products by $L$ only require each node to communicate with its neighbours, the same is not true with $L^{\frac{1}{2}}$, given that $L^\frac{1}{2}$ need not be compatible with the graph topology). However, as shown below, eliminating the second variable yields \eqref{ParametricAlgorithm} and, for the particular choices of $\eta=1,\beta^2=\frac{1}{2}$, and $L=I-W$, Algorithm 1. Consider  two consecutive  $z$-updates, i.e.,
\begin{align*}
&z^{k+2}=z^{k+1}+\alpha R(z^{k+1})+\beta L^{\frac{1}{2}}w^{k+1}-\eta Lz^{k+1}\\
&z^{k+1}=z^{k}+\alpha R(z^{k})+\beta L^{\frac{1}{2}}w^k-\eta Lz^{k}.
\end{align*}
Consider their difference
\begin{align*}
z^{k+2}&=2z^{k+1}-z^k+\beta L^{\frac{1}{2}}(w^{k+1}-w^k)-\eta L(z^{k+1}-z^k)\\
&+\alpha\big(R(z^{k+1})-R(z^k)\big).
\end{align*}
Note that $w^{k+1}-w^k=-\beta L^{\frac{1}{2}}z^k$ and observe that  its substitution in the $z^{k+2}$-update results in the elimination of $w^{k+1}$ and $w^k$ and the disappearance of $L^{\frac{1}{2}}$, i.e.,
\begin{align*}
z^{k+2}&=(2I-\eta L)z^{k+1}-(I+\beta
^2L-\eta L)z^k\\
&+\alpha\big(R(z^{k+1})-R(z^k)\big),
\end{align*}
which is precisely \eqref{ParametricAlgorithm}. We remark that a similar idea of using $L^\frac{1}{2}$ has been employed in \cite{uribe2020dual} and, for further insights, see section 3 therein.

\subsection{Convergence Theorems}\label{SectionTheoremStatements}
The following two theorems establish the convergence results. The first is a result about the map $\tilde{F}$ and the second about the map $F$.

\ 

\begin{theorem}\label{ImportantTheorem} Let $\lambda_1, \ldots, \lambda_{dN}$ be the eigenvalues of $L$ and choose $\eta>0$ and $\beta>0$ such that the two complex roots $\gamma_1(\lambda_s,\eta,\beta)$ and $\gamma_2(\lambda_s,\eta,\beta)$ of all the $dN$ polynomial equations
\begin{align*}
x^2+\lambda_s\eta x+ \lambda_s\beta^2=0
\end{align*}
satisfy the two following conditions,
\begin{description}
\item[1)] $|1+ \gamma_j(\lambda_s,\eta,\beta)|\leq 1$ 
\item[2)] $|1+\gamma_j(\lambda_s,\eta,\beta)|=1$, if and only if $\lambda_s=0$,
\end{description}
for any $s=1, \ldots, dN$ and $j=1,2$. Then, there exists $\alpha^\star$ such that, for $0<\alpha<\alpha^\star$,
\begin{align*}
\rho\big(\mathbf{J}_{\tilde{F}}(\psi(x^\star)\big)<1.
\end{align*}
\end{theorem}

\begin{rem}
 Note that the choices $\eta=1$ and $\beta^2=\frac{1}{2}$ (which lead to Algorithm 1) satisfy the conditions of Theorem 1. In fact, since  $0\leq \lambda_{s}<2$, 
\begin{align*}
\gamma_{j}(\lambda_{s},1,\frac{1}{2})\in \Big\lbrace \frac{-\lambda_{s} \pm i\sqrt{\lambda_{s}}\sqrt{2-\lambda_{s}}}{2}\Big\rbrace
\end{align*}
for $j=1,2$, and, hence,
\begin{align*}
|1+\gamma_{j}(\lambda_{s},\eta,\beta)|^2&=\bigl(1-\frac{\lambda_{s}}{2}\bigr)^2+\frac{2\lambda_{s}-\lambda_{s}^2}{4}\\
&=1-\frac{\lambda_{s}}{2}.
\end{align*}
The fact that $0\leq \lambda_{s}<2$ is enough to yield conditions 1) and 2) of the theorem.
\end{rem}

\

In order to state the next theorem, consider the  matrix  
\begin{align}
\hat{U}=\frac{1}{\sqrt{N}}\mathbf{1}_N\otimes I_d \label{hatU},
\end{align}
the columns of which form an orthonormal basis of $\mbox{ker}(L)$.

\ 

\begin{theorem} \label{Corolary}
Let $\alpha, \eta$, and $\beta$ be such that $\rho\big(\mathbf{J}_{\tilde{F}}(\psi(x^\star)\big)<1$ (the existence of such a choice is ensured by Theorem \ref{ImportantTheorem}). Let $(z^{k},w^{k})$, for $k=0,1,...$,  be a trajectory of the iteration \eqref{Picard2},  initialized at $(z^0,w^0)\in \mathbb{R}^{dN}\times \mathbb{R}^{dN}$. Then, there exists a neighborhood $\mathcal{V}$ of $\psi(x^\star)$ such that if $(z^k,\tilde{U}^Tw^k) \in \mathcal{V}$ for some $k$, then this trajectory converges to
\begin{align*}
\begin{bmatrix}
\mathbf{1}_N \otimes x^\star \\
\hat{U}\hat{U}^Tw^0-\frac{\alpha}{\beta}(L^\frac{1}{2})^+R(\mathbf{1}_N \otimes x^\star)
\end{bmatrix},
\end{align*}
with at least linear rate.
\end{theorem}

\subsection{Proof of Theorem \ref{ImportantTheorem}} \label{Proof1}
We will show that we can trap the eigenvalues of $\mathbf{J}_{\tilde{F}}(\psi(x^\star))$ in the open ball in $\mathbb{C}$ of center $0$ and radius $1$ by  adjusting $\alpha, \beta,$ and $\eta$. The Jacobian of $\tilde{F}$ at $\psi(x^\star)$ is given by 
\begin{align*}
\mathbf{J}_{\tilde{F}}\big(\psi(x^\star)\big)=I+\hat{A}(\eta,\beta)+B(\alpha),
\end{align*}
 where
\begin{align*}
\hat{A}(\eta,\beta)&=\begin{bmatrix}
-\eta L & \beta L^\frac{1}{2}\tilde{U}\\
-\beta\tilde{U}^TL^\frac{1}{2} & \mathbf{0}
\end{bmatrix},
\end{align*}
and
\begin{align*}
B(\alpha)&=\begin{bmatrix}
\alpha \mathbf{J}_R(\mathbf{1}_N\otimes x^\star)&\mathbf{0}\\
\mathbf{0}&\mathbf{0}
\end{bmatrix}.
\end{align*}

We divide the proof in two lemmas. Lemma 1 shows that a choice of $\eta>0$ and $\beta>0$ according to the conditions of the Theorem implies that the matrix $I+\hat{A}(\eta,\beta)$ has  a semisimple eigenvalue equal to $1$, and that the remaining eigenvalues have magnitude less than $1$. Lemma 2  deals with $\alpha$ and shows, by appealing to a result on the perturbation of semisimple eigenvalues, that this parameter can be tuned to yield the conclusion of  Theorem 1.

\

\begin{lemma} 
Let $\eta>0$ and $\beta>0$ be chosen according to the conditions of Theorem 1. Then 
\begin{description}
\item[i)] $1$ is a semisimple eigenvalue of $I+\hat{A}(\eta,\beta)$ with multiplicity $d$ and corresponding eigenspace given by $\ker(L)\times \lbrace 0\rbrace$;
\item[ii)] All the remaining eigenvalues have magnitude less than $1$. 
\end{description}
\end{lemma}

\begin{proof}

 We start by analyzing the eigenvalues of $\hat{A}(\eta,\beta)$; towards that goal, consider the eigenvalues of
\begin{align*}
A(\eta,\beta)=\begin{bmatrix}
-\eta L & \beta L^\frac{1}{2} U\\
-\beta UL^\frac{1}{2} & \mathbf{0}
\end{bmatrix},
\end{align*}
where $U=[\tilde{U},\hat{U}]$ ($\hat{U}$ was defined in \eqref{hatU} and $\tilde{U}$ in the first paragraph of Section \ref{Fmaps}). Since
\begin{align*}
A(\eta,\beta)=\begin{bmatrix}
\hat{A}(\eta,\beta)&\mathbf{0}\\
\mathbf{0}&\mathbf{0}
\end{bmatrix},
\end{align*}
the non-zero eigenvalues of $\hat{A}(\eta,\beta)$ coincide with the non-zero eigenvalues of $A(\eta,\beta)$.

Let $V$ be an orthogonal matrix such that $V^TLV=\Lambda$, where $\Lambda=\mbox{diag}(\lambda_1, \ldots, \lambda_{dN})$,  and consider the following unitary similarity
\begin{align}\label{TensorPoductMatrix}
\begin{bmatrix}
V^T&\mathbf{0}\\
\mathbf{0}&V^TU
\end{bmatrix}A(\eta,\beta)\begin{bmatrix}
V&\mathbf{0}\\
\mathbf{0}&U^TV
\end{bmatrix}=\begin{bmatrix}
-\eta \Lambda& \beta\Lambda^\frac{1}{2} \\
-\beta \Lambda^\frac{1}{2} &\mathbf{0}
\end{bmatrix}.
\end{align}
The eigenvalues are preserved by unitary similarity and, hence, the eigenvalues of \eqref{TensorPoductMatrix} coincide with those of $A(\eta,\beta)$.

 Let $\xi$ be an eigenvalue of \eqref{TensorPoductMatrix} and $(u,v)\neq (0,0)$ be an associated eigenvector. There must exist an $s$ such that  $u_s\neq 0$ or $v_s\neq 0$, which implies that $\xi$ is an eigenvalue of the $2\times 2$ matrix
\begin{align*}
\begin{bmatrix}
-\eta \lambda_{s}&\beta \sqrt{\lambda_{s}}\; \\
-\beta\sqrt{\lambda_{s}}&0
\end{bmatrix},
\end{align*}
thus $\xi$ is a solution of the second-degree equation
\begin{align}\label{ChacteristicPoly2}
x^2+\eta \lambda_{s} x+\lambda_{s}\beta^2=0.
\end{align}
Note that the converse also holds, meaning that any root of \eqref{ChacteristicPoly2} gives rise to an eigenvalue of \eqref{TensorPoductMatrix}. This completely characterizes the eigenvalues of \eqref{TensorPoductMatrix}.

Now if $\eta>0$ and $\beta>0$  are chosen according to the conditions of Theorem 1, then, from the derivation above, 

\begin{align*}
\rho\big(I+\hat{A}(\eta,\beta)\big)\leq 1.
\end{align*}
Moreover, if $\lambda$ is an eigenvalue of $I+\hat{A}(\eta,\beta)$ associated to the eigenvector $(u,v)$, such that $|\lambda|=1$, it follows, from condition 2) of  Theorem 1, that $\lambda=1$ and, therefore, $(u,v) \in \ker (\hat{A}(\eta,\beta))$. Additionally, observe that: 
\begin{itemize}
	\item[a)]  $\ker (\hat{A}(\eta,\beta))=\ker (L) \times \lbrace 0 \rbrace$ (the columns of $\tilde{U}$ form a basis for  $\mbox{range}(L)$);
	\item[b)] the orthogonal of $\ker (\hat{A}(\eta,\beta))$ is equal to $\mbox{range}(L) \times \mathbb{R}^{d(N-1)}$; 
	\item[c)] both these spaces are invariant under $\hat{A}(\eta,\beta)$. 
\end{itemize}	

From these observations, we conclude that $1$ is a semisimple eigenvalue of $I+\hat{A}(\eta,\beta)$ with multiplicity $d$ and corresponding eigenspace given by $\ker (L) \times \lbrace 0 \rbrace$, and that
 all the other eigenvalues have magnitude less than $1$.
\end{proof}

\

\begin{lemma} Let $\eta>0$ and $\beta>0$ be chosen according to the conditions of Theorem 1. Then, there exists $\alpha^\star$ such that, for $0<\alpha<\alpha^\star$, 
\begin{align*}
\rho\big(\mathbf{J}_{\tilde{F}}(\psi(x^\star)\big)<1.
\end{align*}
\end{lemma}
\begin{proof}
Let $C(\alpha)=I+\hat{A}(\eta,\beta)+B(\alpha)$, which can be seen as a continuous perturbation of the matrix $C(0)=I+\hat{A}(\eta,\beta)$ (recall that $B(0) = \mathbf{0}$). We will show that choosing a sufficiently small $\alpha$ traps the eigenvalues of $C(\alpha)$ inside the ball in $\mathbb{C}$ of center $0$ and radius $1$. 

Theorem \ref{ContinuousEigenvaluesCurve} (included in Appendix \ref{Auxiliary results Section}, for convenience; see also \cite{kato2013perturbation}) implies that there are $2dN-d$ continuous and complex-valued functions $\lambda_1(\alpha),\ldots, \lambda_{2dN-d}(\alpha)$, such that $\lbrace \lambda_s(\alpha) : s=1, \ldots, 2dN-d\rbrace$ is the set of eigenvalues of $C(\alpha)$. By reordering, if necessary, we may assume that, for $s=1, \ldots, d$, $\lambda_s(0)=1$ (recall, from Lemma 1, that $1$ is a semisimple eigenvalue of $C(0)$). With this choice of order, it follows that  $|\lambda_s(0)|<1$, for $d+1\leq s\leq 2dN-d$.

By continuity, there exists $\alpha_1>0$ ensuring that $|\lambda_s(\alpha)|<1$, for all $0\leq \alpha< \alpha_1$ and for $d+1\leq s\leq 2dN-d$.

Theorem \ref{DifferentialEigenvaluesCurve}  (included in Appendix \ref{Auxiliary results Section}, for convenience; see also \cite{lancaster1964eigenvalues})  implies that $\lambda_1(\alpha), \ldots, \lambda_{d}(\alpha)$ are, not only continuous, but also differentiable at $0$. Moreover, their derivatives at $0$ are among the eigenvalues of 
\begin{align*}
\begin{bmatrix}
\hat{U}^T &\mathbf{0} \end{bmatrix}\biggl(\frac{dC}{d\alpha}\Bigr|_{\alpha=0}\biggr)\begin{bmatrix}
\hat{U}^T\\
\mathbf{0}
\end{bmatrix}=\hat{U}^T\mathbf{J}_R(\mathbf{1}_N\otimes x^\star)\hat{U}.
\end{align*}
We observe that $\hat{U}^T\mathbf{J}_R(\mathbf{1}_N\otimes x^\star)\hat{U}=\mathbf{J}_H(x^\star)-I$.

Condition \eqref{JacobianCondition}, i.e., $\rho\big(\mathbf{J}_H(x^\star)\big)<1$, has not yet entered the game and it is here that it plays a crucial role. In fact, observe that condition \eqref{JacobianCondition} implies that the eigenvalues of $\mathbf{J}_H(x^\star)-I$ have negative real part.

 A geometrical argument now settles the question: each $\lambda_s(\alpha)$, for $s=1, \ldots, d$, can be viewed as a continuous curve in $\mathbb{R}^2$ and can, therefore, be written as $\lambda_s(\alpha)=(\lambda_{s}^{(1)}(\alpha),\lambda_{s}^{(2)}(\alpha))$. So far, we proved that these curves satisfy
\begin{align*}
\lambda_s(0)&=\begin{bmatrix} 1\\0\end{bmatrix}
\end{align*}
and
\begin{align*}
\frac{d\lambda_{s}^{(1)}}{d\alpha}\Big|_{\alpha=0}&<0.
\end{align*}
By observing that
\begin{align*}
\frac{d\, \|\lambda_s(\alpha)\|^2}{d\alpha}\Big|_{\alpha=0}&=2\, \frac{d\lambda_{s}^{(1)}}{d\alpha}\Big|_{\alpha=0}<0 
\end{align*}
and
\begin{align*}
\|\lambda_s(0)\|^2&=1,
\end{align*}
we conclude that there must exist an $\alpha_2>0$ such that 
\begin{align*}
\|\lambda_s(\alpha)\|^2<1,
\end{align*}
for $s=1, \ldots, d$ and for $0<\alpha<\alpha_2$. Finally, choosing $\alpha^\star=\min \lbrace \alpha_1, \alpha_2\rbrace$, it follows that, for $0<\alpha<\alpha^\star$,
\begin{align*}
\rho(C(\alpha))<1. &&& \qedhere
\end{align*} 
\end{proof}
\color{black}

\subsection{Proof of Theorem \ref{Corolary}}  \label{Proof2}
For any vector $w \in \mathbb{R}^{dN}$,  let $\tilde{w}:=\tilde{U}\tilde{U}^Tw$ and $\bar{w}:=\hat{U}\hat{U}^Tw$. Observe that $w=\bar{w}+\tilde{w}$ is the orthogonal decomposition of $w$ with components in $\ker (L)$ and $\mbox{range}(L)$.

Consider the sequence $(z^k,w^k)$, $k=0,1,...$, produced by the iteration \eqref{Picard2}, with  a given initialization $(z^0,w^0)$. We note that\footnote{There is a slight abuse of notation here: we are using $\bar{w}^k$ and $\tilde{w}^k$ as, respectively, $\hat{U}\hat{U}^Tw^k$ and $\tilde{U}\tilde{U}^Tw^k$.} $\bar{w}^k=\bar{w}^0$, for all $k$, and $z^k$ only depends on $\tilde{w}^k$. From these two facts, it follows that the sequence produced by the iteration
\begin{align}
\begin{bmatrix}
u^{k+1}\\
v^{k+1}
\end{bmatrix}=F\Big(\begin{bmatrix}
u^k\\
v^k
\end{bmatrix}\Big)\ \label{Picard3}
\end{align}
with initialization $u^0=z^0$ and $v^0=\tilde{w}^0$, must satisfy
\begin{align*}
z^k&=u^k \\
w^k&=\bar{w}^0+v^k,
\end{align*}
for all $k$. Moreover, for all $k$, $v^k=\tilde{v}^k$.

Said differently, the observations above show that to understand the trajectories followed by $(z^{k},w^{k})$, given an arbitrary initialization $(z^0,w^0)$, we can just study the trajectories initialized in $\mathbb{R}^{dN}\times \mbox{range}(L)$. Moreover, the second component of a trajectory initialized in $\mathbb{R}^{dN}\times \mbox{range}(L)$ stays in $\mbox{range}(L)$. 

Furthermore, $\tilde{U}$ establishes an isomorphism between $\mathbb{R}^{dN}$ and $\mbox{range}(L)$, with inverse $\tilde{U}^T$. With this in mind, let $(z^0, w^0)\in \mathbb{R}^{dN}\times \mbox{range}(L)$ and observe that
\begin{align*}
z^{k+1}&=z^k+\alpha R(z^k)+\beta L^\frac{1}{2}\tilde{U}\tilde{U}^Tw^k-\eta Lz^k\\
\tilde{U}\tilde{U}^Tw^{k+1}&=\tilde{U}\tilde{U}^Tw^k-\beta L^\frac{1}{2}z^k.	
\end{align*}
Defining $\tilde{w}^k=\tilde{U}^Tw^k$ we obtain 
\begin{align*}
z^{k+1}&=z^k+\alpha R(z^k)+\beta L^\frac{1}{2}\tilde{U}\tilde{w}^k-\eta Lz^k\\
\tilde{w}^{k+1}&=\tilde{w}^k-\beta \tilde{U}^TL^\frac{1}{2} z^k.	
\end{align*}
We conclude that
\begin{align*}
\begin{bmatrix}
z^k \\ w^k
\end{bmatrix}& = \begin{bmatrix} z^k \\ \tilde{U}\tilde{w}^k  \end{bmatrix}
\end{align*}
and
\begin{align*}
\begin{bmatrix} z^{k+1} \\ \tilde{w}^{k+1}  \end{bmatrix} &= \tilde{F} \Bigl( \begin{bmatrix} z^k \\ \tilde{w}^k  \end{bmatrix} \Bigr).
\end{align*}
These observations, combined with Theorem \ref{ImportantTheorem} and Ostrowski's Theorem (see Appendix \ref{Auxiliary results Section}) are enough to conclude the proof of Theorem \ref{Corolary}.

\section{Intuition for \eqref{ParametricAlgorithm}} \label{Intuition}
To provide some intuitive insight into the rationale that leads to the Banach-Picard iteration  \eqref{Picard2}, and, after the elimination of the second variable, to \eqref{ParametricAlgorithm}, suppose for a moment that the maps $H_n$ in \eqref{DistributedMap} had the form $H_n=I -  \nabla f_n$, for some  functions $f_n : \mathbb{R}^d \to \mathbb{R}$, and note that, in this case, the Banach-Picard iteration with $H $ reduces to the familiar gradient descent with unit step-size for finding a stationary point of $$
f=\frac{1}{N}\sum_{i=1}^N f_n.
$$
Moreover, condition \eqref{JacobianCondition}, in this case, reads
\begin{align*}
\rho\big(\mathbf{J}_H(x^\star)\big) = \rho\big(I-\mathbf{J}^2_f(x^\star)\big) < 1,
\end{align*}
which implies that $\mathbf{J}^2_f(x^\star)\succ 0$, showing that, in this case, $f$ is locally strictly convex, thus $x^\star$ is a local minimum of $f$. 

A common approach to pursue a distributed algorithm to maximize $f$ is to formulate the problem as (see, e.g., \cite{jakovetic2020primal}, \cite{uribe2020dual}, and \cite{mansoori2019general})
\begin{equation}\label{intuition1}
\begin{aligned}
& \underset{z\in\mathbb{R}^{dN}}{\text{minimize}}
& &\sum_{n=1}^N f_n(z_n) \\
& \text{subject to}
& &  \beta' L^\frac{1}{2}z=0
\end{aligned},
\end{equation}
where $z = [z_1^T,..., z_N^T]^T$, $L$ is a matrix satisfying the conditions given in Section \ref{KKT section}), and $\beta' > 0$; the constraint $L^\frac{1}{2}z=0$ imposes that $z_1 = \cdots = z_N$ (as shown in \cite{Mokhtari2016} and  \cite{uribe2020dual}). The so-called augmented Lagrangian  for this problem is 
\begin{align}
\mathcal{L}(z,w) = \sum_{n=1}^N f_n(z_n) + \beta' w^T L^{\frac{1}{2}}z + \frac{\eta'}{2}\|L^{\frac{1}{2}}z\|^2,
\end{align}
where $w$ is the vector or Lagrange multipliers and $\eta' > 0$ (see \cite{bertsekas2014constrained}). The \textit{augmented Lagrangian method}  consists in alternating between minimizing $\mathcal{L}(z,w)$ with respect to $z$ and taking a gradient ascent step with respect to $w$, while keeping the other variable fixed. If instead of exact minimization, the $z$-step is itself a gradient  ascent step, we have the Arrow-Hurwitz-Uzawa method (see \cite{jakovetic2020primal})
\begin{align*}
z^{k+1}=z^k&-\alpha \nabla_z \, \mathcal{L}(z^k,w^k)\\
w^{k+1}=w^k&+\alpha \nabla_w \, \mathcal{L}(z^k,w^k) = w^k + \alpha\beta' L^{\frac{1}{2}}z^k.
\end{align*}
Setting $\alpha \eta^\prime=\eta$ and $\alpha \beta^\prime=\beta$ recovers \eqref{ParametricAlgorithm} for the case where $H_n=I -  \nabla f_n$. What this shows is that \eqref{ParametricAlgorithm} applied to this particular case, can be interpreted as the Arrow-Hurwitz-Uzawa method to find a saddle point of the augmented Lagragian. 

As a final observation, note that, if we again let $H_n=I+\nabla f_n$, then Algorithm 1 is nothing but the algorithm in  \cite{shi2015extra} for the concave case.

\section{Conclusion and Future Work} \label{Conclusions}
The article proposed an algorithm for the distributed computation of attractors of average maps. The conditions of the average map considered were rather minimal, yet sufficiently strong to encapsulate relevant algorithms like the EM algorithm or the algorithm in \cite{sanger1989optimal} for PCA, as it is demonstrated in detail and experimentally assessed in an upcoming article \cite{FranciscodistributedEMPCA}. Regardless of the minimality of the set of conditions, we were, nevertheless, able to provide local guarantees of linear convergence by employing a proof technique based on perturbation theory for linear operators, which, to the best of our knowledge, had not been used in distributed optimization. In an upcoming article \cite{FranciscodistributedEMPCA}, we verify the attractor condition for the EM algorithm and the algorithm in \cite{sanger1989optimal} for PCA, and, therefore, we propose a distributed EM algorithm and a distributed algorithm for PCA.


%


\appendices
\section{Auxiliary Results} \label{Auxiliary results Section}
This appendix contains three theorems that are instrumental in the proofs presented in Sections \ref{Proof1} and \ref{Proof2}.

\ 

\begin{theorem}[Ostrowski's Theorem; for a proof, see \cite{ortega2000iterative}] \label{Ostrowski's THeorem}
Suppose $G:\mathcal{D}\subseteq \mathbb{R}^d \to \mathbb{R}^d$ has a fixed point $x^\star$ in the interior of $D$. Suppose as well that $G$ is $C^1$ in a neighborhood of $x^\star$ and that 
\begin{align*}
\rho\big(\mathbf{J}_G(x^\star)\big)<1.
\end{align*}
Then, there exists a norm $\|\cdot\|$, a ball $\mathcal{B}=\lbrace x \in \mathbb{R}^d: \|x-x^\star\|< \delta\rbrace$, and a number $0<\sigma<1$ such that, for all $x \in \mathcal{B}$,
\begin{align*}
\|G(x)-x^\star\|\leq \sigma\|x-x^\star\|.
\end{align*}

\end{theorem}

\ 

\begin{theorem}[for a proof, see \cite{kato2013perturbation}]\label{ContinuousEigenvaluesCurve}
Let $G(x)$ be an unordered $N$-tuple of complex numbers depending continuously on a real variable $x$ in a (closed or open) interval $I$. Then, there exists $N$ single-valued continuous functions, $\mu_n(x)$, $n=1, \ldots, N$, the values of which constitute the $N$-tuple $G(x)$, for each $x \in I$.
\end{theorem}

\ 

\begin{theorem}[for a proof, see \cite{lancaster1964eigenvalues}]\label{DifferentialEigenvaluesCurve}
Let $A$ and $B$ be $n\times n$ square matrices and consider the parametric family $C(\alpha)=A+\alpha B$. Let $\lambda$ be a semisimple eigenvalue of $A$ with multiplicity $p$ (geometric and algebraic multiplicities coincide) and suppose that there is an orthogonal $n\times n$ matrix $U=[\hat{U},\tilde{U}]$, where $\hat{U}$ and $\tilde{U}$ have, respectively, dimensions $n\times p$ and $n\times (n-p)$, such that
\begin{align*}
U^TAU=\begin{bmatrix}
\lambda I_p &\mathbf{0}\\
\mathbf{0} & \tilde{A}
\end{bmatrix},
\end{align*}
and where $\lambda$ is not an eigenvalue of $\tilde{A}$. Let $\mu_1(\alpha),\ldots,\mu_p(\alpha)$ be the $p$ eigenvalue curves of $C(\alpha)$ that satisfy $\mu_j(0)=\lambda$, for $j=1,\ldots, p$ (the existence of these curves is guaranteed by Theorem \ref{ContinuousEigenvaluesCurve}). Then, each $\mu_j$ is differentiable at $0$ and the derivatives of $\mu_1(\alpha),\ldots,\mu_p(\alpha)$ at $0$ are among the eigenvalues of $\hat{U}^TB\hat{U}.$
\end{theorem}

\section*{Acknowledgment}
This work was partially funded by the Portuguese \textit{Fundação para a Ciência e Tecnologia} (FCT), under grants  PD/BD/135185/2017 and UIDB/50008/2020.  João Xavier was supported in part by the Fundação para a Ciência e Tecnologia, Portugal, through the Project LARSyS, under Project FCT Project UIDB/50009/2020 and Project HARMONY PTDC/EEI-AUT/31411/2017 (funded by Portugal 2020 through FCT, Portugal, under Contract AAC n 2/SAICT/2017–031411. IST-ID funded by POR Lisboa under Grant LISBOA-01-0145-FEDER-031411).

\ifCLASSOPTIONcaptionsoff
  \newpage
\fi

\end{document}